\documentclass[10pt]{article}
\usepackage{hyperref}
\usepackage{fullpage,graphicx,graphics,amsmath,amsthm}
\usepackage{epstopdf}

\newtheorem{theorem}{Theorem}[section]
\newtheorem{lemma}[theorem]{Lemma}

\numberwithin{figure}{section}
\newcommand{\del}{\backslash}

\begin{document}

\title{Strengthened chain theorems for different versions of 4-connectivity}

\author{
Guoli Ding \\
Department of Mathematics, Louisiana State University, Baton Rouge, USA
\and
Chengfu Qin \thanks{Corresponding author}\\
Department of Mathematics, Nanning Normal University, Nanning, China}

\date{\today}

\maketitle

\begin{abstract}
The chain theorem of Tutte states that every 3-connected graph can be constructed from a wheel $W_n$ by repeatedly adding edges and splitting vertices. It is not difficult to prove the following strengthening of this theorem: every non-wheel 3-connected graph can be constructed from $W_4$ by repeatedly adding edges and splitting vertices. In this paper we similarly strengthen several chain theorems for various versions of 4-connectivity.\end{abstract}

\section{Introduction}

If a graph $G$ contains a minor $H$ then by a {\it $(G,H)$-chain} we mean a sequence $G_0,G_1,...,G_t$ of graphs such that $G_0=G$, $G_t=H$, and each $G_i$ $(1\le i\le t)$ is a minor of $G_{i-1}$. If $\cal H$ is a class of graphs then by a {\it $(G,\cal H)$-chain} we mean a $(G,H)$-chain for some $H\in\cal H$. A result that asserts the existence of a chain of certain type is often called a {\it chain theorem}. Probably the best known chain theorem is the following result of Tutte. Let $\mathcal W=\{W_n:n\ge3\}$, where $W_n$ is the wheel on $n+1$ vertices. Let $|G|$ and $\|G\|$ denote the number of vertices and edges, respectively, of a graph $G$.

\begin{theorem}[\cite{tutte}] \label{thm:tutte}
For every $3$-connected graph $G$ there exists a $(G,\cal W)$-chain $G_0,G_1,...,G_t$ of $3$-connected graphs such that $\|G_i\|=\|G_{i-1}\|-1$ for all $i=1,...,t$.
\end{theorem}

This theorem is a very powerful induction tool for analyzing 3-connected graphs since it says that every 3-connected graph can be reduced to a wheel, by reducing only one edge at a time, while maintaining 3-connectivity. If we reverse this reduction process we see that every 3-connected graph can be constructed from a wheel by repeatedly performing ``undeletion" and ``uncontraction" (note that both operations preserve 3-connectivity). When such a chain theorem is used to generate all 3-connected graphs, it is desirable to limit the number of starting graphs (wheels). This in fact can be done. It is not difficult to show that every non-wheel 3-connected graph can be constructed using the same operations starting from $W_4$, instead of a general unknown wheel. For completeness we include a proof of this result at the end of the paper. From the algorithmic point of view, this strengthened chain theorem is much more efficient since, in order to generate all 3-connected graphs, there is no need for us to consider extensions of all wheels, it is enough for us to only consider extensions of a single wheel $W_4$.

Currently, there are no known chain theorems for 5- or higher connectivity. However, there are several results for different versions of 4-connectivity. The purpose of this paper is to obtain similar strengthenings for three chain theorems, respectively, for 4-connectivity, weakly 4-connectivity, and quasi 4-connectivity. In some cases we also refine the involved operations.

In the next three sections we discuss these three connectivities separately. Our main results are Theorem \ref{thm:main}, Theorem \ref{thm:w4c}, and Theorem \ref{thm:q4c+}. Before closing this section we introduce a few definitions.  Let $G=(V,E)$ be a graph. For any $X\subseteq V$, let $N_G(X)=\{y\in V\del X$: there exists $x\in X$ with $xy\in E\}$. Suppose $x\in V$. We will write $N_G(x)$ for $N_G(\{x\})$. As usual, $|N_G(x)|$ is the {\it degree} of $x$, which is denoted by $d_G(x)$. We will drop the subscript $G$ if there is no confusion.

A {\it separation} of a graph $G=(V,E)$ is a pair of subgraphs $G_i=(V_i,E_i)$ ($i=1,2$) of $G$ such that $G_1\cup G_2=G$, $E_1\cap E_2=\emptyset$, and $V_i\ne V$ for both $i=1,2$. If $|V_1\cap V_2|=k$ then the separation is called a {\it $k$-separation} and $V_1\cap V_2$ is called a {\it $k$-cut}. For any integer $k\ge3$, note that $G$ is $k$-connected if and only if $G$ is simple, $|G|>k$, and $G$ has no $k'$-separation for any $k'<k$.

For any two subgraphs $G_1,G_2$ of a graph $G$, we will write $G_1\del G_2$ for $G_1\del V(G_2)$. We will use $si(G)$ to denote the simplification of $G$. If $e$ is an edge of $G$, we use $\bar e$ to denote the vertex of $G/e$ obtained from contracting $e$.
Further, if $e$ is contained in exactly one triangle which contains three edges $\{e,f,h\}$, then we see that $si(G/e)\cong G/e\backslash f$. In the rest of the paper, $si(G/e)$  will be replaced by $G/e\backslash f$ if $e$ is contained in exactly one triangle and  $si(G/e)$  will be replaced by $G/e$ if $e$ is contained in no triangle.
See \cite{diestel} for undefined notation and terminology.

\section{On 4-connected graphs}

Let $\mathcal C = \{C^2_n: n\ge5\}$, where $C^2_n$ is the graph obtained from the $n$-cycle $C_n$ by joining vertices of distance two on the cycle. Note that each $C^2_n$ is 4-connected, and in addition, $C^2_5=K_5$ and $C^2_6=K_{2,2,2}=$ octahedron. A cubic graph with at least six vertices is called {\it internally $4$-connected} if its line graph is 4-connected. Let $\mathcal L=\{L: L$ is the line graph of an internally 4-connected cubic graph\}. The following is a well known chain theorem for 4-connected graphs.

\begin{theorem}[\cite{fontet} and \cite{martinov}] \label{thm:11}
For every $4$-connected graph $G$ there exists a $(G, \cal C\cup L)$-chain $G_0,G_1,...,G_t$ of $4$-connected graphs such that $G_i=si(G_{i-1}/e_i)$ for all $i=1,...,t$.
\end{theorem}

This chain theorem has been strengthened as follows.

\begin{theorem}[\cite{qin}] \label{thm:qin}
For every $4$-connected graph $G\not\in\cal C\cup L$ there exists a $(G,\{C_5^2,C_6^2\})$-chain $G_0,G_1,...,G_t$ of $4$-connected graphs such that $G_i=si(G_{i-1}/e_i)$ for all $i=1,...,t$. In addition, if $G$ is nonplanar then $G_t=C_5^2$.
\end{theorem}

This result implies that every 4-connected graph $G\not\in {\cal L}\cup {\cal L}$ must contain a 4-connected minor on $|G|-1$ vertices. With this formulation in mind, we can say that  both Theorem \ref{thm:11} and Theorem \ref{thm:qin} are extension of Tutte Theorem to 4-connected graphs. However, this extension is not totally satisfactory since in these two theorems the gap $\|G_{i-1}\|-\|G_{i}\|$ could be arbitrarily large while the gap in Tutte Theorem is only 1. Our main result in this section is to bound this gap.

\begin{figure}[htb]
\centerline{\includegraphics[scale=0.35]{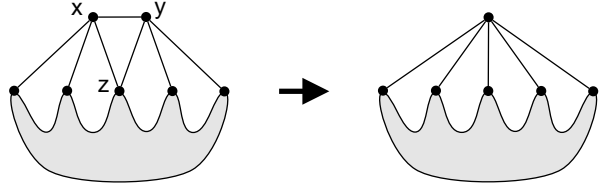}}
\caption{Operation $O_1$, where $d(x)=d(y)=4$.}
\label{fig:split}
\end{figure}

For any integer $n\ge3$, let $B_n$ ({\it biwheel}) be the graph obtained from $C_n$ by adding two nonadjacent vertices ({\it hubs}) and joining them to all vertices of $C_n$ ({\it rim}). Let $B_n^+$ be obtained from $B_n$ by adding an edge ({\it axle}) between the two hubs. Note that $B_3^+=K_5=C_5^2$, $B_4=C_6^2$, and $B_4^+=K_{1,1,2,2}$. In addition, every $B_n$ is planar while every $B_n^+$ is nonplanar. Let $\mathcal B=\{B_n, B_n^+:n\ge4\}$. Then every $B\in {\cal B}$ is 4-connected. Finally, let $O_1$ be the operation illustrated in Figure \ref{fig:split}. That is, $O_1(G)=G/xy\del xz$, where $d(x)=d(y)=4$, and $xyz$ is the only triangle containing $xy$. The following is the main theorem of this section.

\begin{theorem} \label{thm:main}
For every $4$-connected graph $G\not\in\cal B\cup C\cup L$ there exists a $(G,\{B_4^+,B_5\})$-chain $G_0,G_1,...,G_t$ of $4$-connected graphs such that, for each $i=1,...,t$, either $\|G_i\|=\|G_{i-1}\|-1$ or $G_i=O_1(G_{i-1})$. In addition, if $G$ is nonplanar then $G_t=B_4^+$.
\end{theorem}

The following is a useful lemma, which says that ``splitting" a vertex in a $k$-connected graph also results in a $k$-connected graph.
This lemma suggests that Theorem \ref{thm:main} can be formulated in a slightly stronger form. Together with this lemma, Theorem \ref{thm:main} implies that a graph $G\not\in\cal B\cup C\cup L$ is 4-connected if and only if there exists a $(G,\{B_4^+,B_5\})$-chain $G_0,G_1,...,G_t$ such that $\delta(G_i)\ge4$ ($0\le i\le t$) and either $\|G_i\|=\|G_{i-1}\|-1$ or $G_i=O_1(G_{i-1})$ ($1\le i\le t$). This formulation is very useful for generating 4-connected graphs since it ensures that graphs generated according to the theorem are precisely 4-connected graphs that are not contained in $\cal B\cup C\cup L$.

\begin{lemma}\label{lem:split}
Suppose a simple graph $G$ has an edge $xy$ such that $\min\{d_G(x),d_G(y)\}\ge k$ and $si(G/xy)$ is $k$-connected. Then $G$ is also $k$-connected.
\end{lemma}

\noindent Proof. Suppose on the contrary that $G$ is not $k$-connected. Then $G$ has a $k'$-separation $(G_1,G_2)$ with $k'<k$. Assume without loss of generality that $xy\in G_1$. If $G_1\cap G_2$ contains both $x$ and $y$ then $G/xy$ admits a $(k'-1)$-cut, contradicting the $k$-connectivity of $si(G/xy)$. Hence, we may assume $x\in G_1\del G_2$. If $G_1\del G_2$ contains a vertex different from $x$ then $(G_1/xy,G_2)$ is a $k'$-separation of $G/xy$, again contradicting the $k$-connectivity of $si(G/xy)$. Consequently, $x$ is the only vertex of $G_1\del G_2$, which is impossible as $d_G(x)\ge k$. This contradiction proves the lemma. \qed

We also need two other small technical lemmas for proving Theorem \ref{thm:main}.

\begin{lemma}\label{lem:deg}
Let $G$ be a $k$-connected graph with a triangle $xyz$. If $si(G/xy)$ is $k$-connected but $G\del xz$ is not, then $d_G(x) = k$ and $d_G(z)>k$.
\end{lemma}

\noindent Proof. First, since $si(G/xy)$ is $k$-connected, we must have $d_G(z)=d_{si(G/xy)}(z)+1\ge k+1$. Next, suppose on the contrary that $d_G(x)>k$. Let us consider a $k'$-separation $(G_1,G_2)$ of $G\del xz$ with $k'<k$. Without loss of generality, let $x\in G_1$. Then $x\in G_1\del G_2$ and $z\in G_2\del  G_1$, because otherwise $(G_1,G_2)$ would induce a $k'$-separation of $G$, which violates the $k$-connectivity of $G$. As a result, $y\in G_1\cap G_2$. Since $d_G(x)>k$, $G_1\del G_2$ contains a vertex different from $x$. It follows that $(G_1,G_2)$ induces a $k'$-separation of $G/xy$, which violates the $k$-connectivity of $si(G/xy)$. This contradiction proves the lemma. \qed

\begin{lemma}\label{lem:side}
Let $G$ and $si(G/xy)$ be $k$-connected. If $N_G(x)\del N_G(y)=\{y,z\}$ then $si(G/xz)$ is $k$-connected.
\end{lemma}

\noindent Proof. Suppose otherwise. Then $G$ has a $k$-separation $(G_1,G_2)$ with $x,z\in G_1\cap G_2$. Since $si(G/xy)$ is $k$-connected, we have $y\not\in G_1\cap G_2$. By symmetry, let $y\in G_1\del G_2$. Then $x$ is not adjacent to any vertex of $G_2\del G_1$, which implies that $V(G_1\cap G_2)\del x$ is a $(k-1)$-cut of $G$, a contradiction. \qed

The next lemma is the main part of our proof of Theorem \ref{thm:main}. We consider the following two assertions.

\noindent\makebox[9mm][r]{($*$)} \ {\it $G$ has a $4$-connected minor $H=G\del  e$, $G/e$, or $O_1(G)$ such that $H\not\in \mathcal C\cup \mathcal L\cup\mathcal B\del\{B_5,B_4^+\}$. \\  \makebox[10mm]{} Moreover, $H$ is not planar if $G$ is not.}\\
\makebox[9mm][r]{($**$)} \ {\it $si(G/e)$ is $4$-connected, and $G/e$ is not planar if $G$ is not.}

\begin{lemma}\label{lem:induction}
Suppose a $4$-connected graph $G$ does not satisfy $(*)$. If $e=xy$ is an edge of $G$ that satisfies $(**)$ then $G\del\{x,y\}$ has four distinct vertices $x'$, $y'$, $z$, $z'$ such that $N_G(x)=\{x',y,z,z'\}$, $N_G(y)=\{x,y',z,z'\}$, and both $xx',yy'$ satisfy $(**)$.
\end{lemma}

\noindent Proof. We divide the proof into a few claims.

\noindent {\bf Claim 1.} {\it $e$ must be contained in a triangle.}

Suppose $e$ belongs to no triangle. It follows that $d_{G/e}(\bar{e})\geq 6$.  Then $G/e$ is not 4-regular and thus $G/e\not\in\cal C\cup L$. If $G/e\not\in\mathcal B\del\{B_5,B_4^+\}$ then, by ($**$), $G$ satisfies ($*$) with $H=G/e$, a contradiction. Hence, we must have $G/e=B\in\mathcal B\del\{B_5,B_4^+\}$. Since  $d_{G/e}(\bar{e})\geq 6$, we see that $B=B_{n+1}$ or $B_n^+$, for some $n\ge5$, and $\bar e$ is a hub of $B$. In other words, $G$ is obtained from $B_{n+1}$ or $B_n^+$ by splitting a hub. Let $uv$ be a rim edge of $B$ so that neither $uvx$ nor $uvy$ is a triangle. Let $H=si(G/uv)$.

We first prove that $H$ is 4-connected. From the choice of $uv$ we have $d_H(x)=d_G(x)$ and  $d_H(y)=d_G(y)$, which imply $\min\{d_H(x),d_H(y)\}\ge4$. In addition, note that $si(H/e)=si(G/e/uv)=si(B/uv)=B_n$ or $B_{n-1}^+$, which are 4-connected. Therefore, by Lemma \ref{lem:split}, $H$ is 4-connected.\vspace{2pt}\\
\indent Next, we observe that, in $G$, both $u$, $v$ have degree 4 and $uv$ is contained in only one triangle $uvh$, where $h$ is the hub of $B$ with $h\neq\bar{e}$. Hence, $H=O_1(G)$. From $d_H(\bar{uv})=5$ we also see that $H\not\in \mathcal C\cup \mathcal L\cup\mathcal B\del\{B_5,B_4^+\}$.

Since $G$ does not satisfy $(*)$, we deduce that $H$ is planar but $G$ is not. However, since $e$ satisfies $(**)$, $G/e$ must be nonplanar, implying that $G/e=B_n^+$. As a result, $H/e = B^+_{n-1}$. Since $n\ge5$, we see that $H/e$, and thus $H$, is in fact nonplanar. This contradiction proves Claim 1.

From now on we assume that $e$ is contained in a triangle $xyz$.

\noindent{\bf Claim 2.} {\it Neither $G\del xz$ nor $G\del yz$ is $4$-connected.}

Suppose otherwise. By symmetry, let us assume that $G'=G\del xz$ is 4-connected. We first prove that $G' \in\mathcal B\del\{B_5,B_4^+\}$. To do so, we verify that $G'$ can almost play the role of $H$ in $(*)$.

\noindent (1) {\it If $G$ is nonplanar then $G'$ is also nonplanar.} \\
\indent To see this, note that $G/e$ is nonplanar, as $e$ satisfies $(**)$. Also note that $si(G/e)=si(G'/e)$, as $xyz$ is a triangle. Consequently, $G'/e$, and thus $G'$, is nonplanar, which proves (1).

\noindent (2) $G'\not\in\cal C\cup L$. \\
\indent Suppose on the contrary that $G'\in\cal C\cup L$. Then every edge of $G'$ belongs to triangle of $G'$. As a result, $G$ contains a triangle $xyz_1$ with $d_G(z_1)=4$. It follows that $d_{si(G/e)}(z_1)=3$, which violates the 4-connectivity assumption on $si(G/e)$. This  contradiction proves (2).

Since $G$ does not satisfy $(*)$, we deduce from (1), (2), and the fact $B_4=C_6^2\in \cal C$
that $G'=B$ for some $B\in\mathcal B\del\{B_4, B_5,B_4^+\}$. In other words, $G=B+xz$, where $B$ contains both $xy$ and $yz$. Since $si(G/e)=si(G\del xz/e)=si(B/e)$ and, by ($**$), $si(G/e)$ is 4-connected, we see that $e$ is not incident with a hub of $B$. That is, $e$ must be a rim edge. Since $xz$ is not an edge of $B$, $z$ can not be a hub and thus $yz$ is also a rim edge. Now we see that $G$ is nonplanar since it is obtained by adding a chord to the rim cycle of $B$. By ($**$), $si(G/e)$ is nonplanar, which means $si(B/e)$ is nonplanar and thus $B=B_n^+$ for some $n\ge5$. Let $f$ be the axle of $B$. Then $H=G\del f$ is 4-connected and nonplanar (as $H=B_n+xz$) and $H\not\in \mathcal C\cup \mathcal L\cup\mathcal B$ (as $H$ contains four vertices of degree exceeding 4). It follows that $H$ satisfies $(*)$, which violates our assumption. This contradiction proves Claim 2.

From Claim 2 and Lemma \ref{lem:deg} we deduce that $d_G(x)=d_G(y)=4$ and $d_G(z)>4$.

\noindent{\bf Claim 3.} {\it $e$ is contained in exactly two triangles.}

If $xy$ belongs to more than two triangles than $|N_G(\{x,y\})|=3$, since $d_G(x)=d_G(y)=4$. It follows that $V(G)=\{x,y\}\cup N_G(\{x,y\})$. This implies $G=K_5$ and thus $si(G/e)$ is not 4-connected, a contradiction.

If $xy$ belongs to only one triangle $xyz$, let $H=si(G/e)$. Then $H=G/xy\del xz$ and thus $H=O_1(G)$. Since $H$ has a vertex $\bar{e}$ of degree five, $H\not\in\mathcal C\cup \mathcal L\cup \mathcal B\del\{B_5,B_4^+\}$. Now, by ($**$), $H$ satisfies ($*$), a contradiction.  Hence,  $e$ is contained in exactly two triangles, which is Claim 3.

Now Claim 3 implies that $G\del\{x,y\}$ has four distinct vertices $x',y',z,z'$ such that $N_G(x)=\{x',y,z,z'\}$ and $N_G(y)=\{x,y',z,z'\}$.

\noindent{\bf Claim 4.} {\it Both $xx'$ and $yy'$ satisfy $(**)$.}

By symmetry, we only need to consider $xx'$. By Lemma \ref{lem:side}, $G'=si(G/xx')$ is 4-connected. If $G'$ is planar then $xyz$ and $xyz'$ are faces of $G'$, as $G'$ is 4-connected. Since $N_G(x)=\{x',y,z,z'\}$ and $x'xy$ is not a triangle of $G$, any multiple edge of $G/xx'$ must be parallel to $xz$ or $xz'$. It follows that the planar embedding of $G'$ can be naturally extended into a planar embedding of $G$, which means $G$ is planar. Hence, we conclude that  $xx'$ satisfies ($**$). This confirms Claim 4 and also completes our proof of the lemma. \qed

\begin{lemma}\label{lem:biwheel}
If a $4$-connected graph $G$ does not satisfy $(*)$ then $G\in\cal B\cup C\cup L$.
\end{lemma}

\noindent Proof. Let us assume $G\not\in\cal C\cup L$, for otherwise we are done. Clearly, we only need to show that $G\in\cal B$. By Theorem \ref{thm:qin}, $G$ has an edge $e=x_1x_2$ satisfying ($**$). Then, by Lemma \ref{lem:induction}, $G\del\{x_1,x_2\}$ has distinct vertices $x_0, y, z, x_3$ such that $N(x_1)=\{x_0,y,z,x_2\}$, $N(x_2)=\{x_1, y, z, x_3\}$, and $x_2x_3$ satisfies ($**$). Next, by applying Lemma \ref{lem:induction} to $x_2x_3$ we deduce that $N(x_3)=\{u,v,x_2,x_4\}$ and $x_3x_4$ satisfies ($**$), where $x_4\not\in N(x_2)\cup\{x_2\}$ and $u,v$ are distinct vertices of $N(x_2)\del\{x_3\}=\{x_1,y,z\}$. To determine $u$ and $v$ we observe from the definition of $x_3$ that $x_3\not\in N(x_1)$. Equivalently, $x_1\not\in N(x_3)$, which implies that $x_1$ is not $u$ or $v$. Therefore, $\{u,v\}=\{y,z\}$ and thus $N(x_3)=\{x_2,y,z,x_4\}$.

In general, suppose we have found in $G$ a path $x_1x_2...,x_k$ ($k\ge4$) and distinct vertices $y,z$ such that $d(x_1)=... = d(x_{k-1})=4$, $\{x_1, ..., x_{k-1}\}\subseteq N(y)\cap N(z)$, and $x_{k-1}x_k$ satisfies $(**)$. By applying Lemma \ref{lem:induction} to $x_{k-1}x_k$ we can either extend this path or conclude that $x_k$ is adjacent to $x_1$, which implies $G=B_k$ or $B_k^+$. This proves the lemma. \qed

\bigskip
\noindent{\bf Proof of Theorem \ref{thm:main}}. Let $G\not\in\cal B\cup C\cup L$ be 4-connected. Then a required chain is obtained by repeatedly applying Lemma \ref{lem:biwheel}. \qed

\section{On weakly 4-connected graphs}

In this paper we also consider two variations of 4-connectivity. Let $G$ be 3-connected. Then $G$ is called {\it weakly $4$-connected} if every 3-separation $(G_1,G_2)$ of $G$ satisfies $\min\{\|G_1\|$, $\|G_2\|\}\le4$. We call $G$ {\it quasi $4$-connected} if every 3-separation $(G_1,G_2)$ of $G$ satisfies $\min\{|G_1|$, $|G_2|\}\le4$. Observe that if $\|G\|\le9$ then $G$ is weakly 4-connected if and only if $G$ is 3-connected, and similarly, if $|G|\le6$ then $G$ is quasi 4-connected if and only if $G$ is 3-connected. It is not difficult to see that weakly 4-connectivity implies quasi 4-connectivity. However, the reverse implication does not hold, as shown by $W_5$. Nevertheless, it is easy to see that a quasi 4-connected graph with ten or more edges is weakly 4-connected if and only if every cubic vertex is contained in at most one triangle. In this section we consider weakly 4-connected graphs. A chain theorem for these graphs has been obtained by Geelen and Zhou \cite{w4c}. To describe this theorem we need some preparations.

\begin{figure}[htb]
\centerline{\includegraphics[scale=0.45]{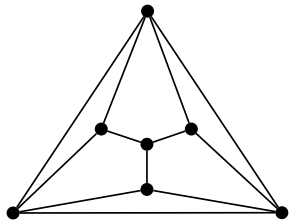}}
\caption{Pyramid $\Pi$.}
\label{fig:pyramid}
\end{figure}

A subgraph of a graph $G$ is called a {\it paw} if it consists of four edges $wx, xy,xz,yz$, and it satisfies $d_G(x)=3$. The graph illustrated in Figure \ref{fig:pyramid}, denoted $\Pi$, is called a {\it pyramid}. It is routine to verify that $\Pi$ is the only weakly 4-connected graph on seven vertices such that its edges can be partitioned into three paws. (One may consider the three possible ways to arrange the three triangles of the paws.)

\begin{figure}[htb]
\centerline{\includegraphics[scale=0.6]{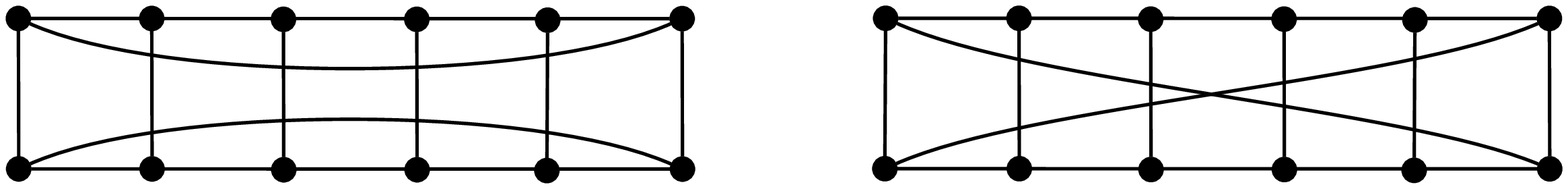}}
\caption{Graphs $A_6$ and $A_6'$.}
\label{fig:a}
\end{figure}

Let $\mathcal B^{\diamond}$ denote the set of planar biwheels $B_n$ for $n\ge4$. For any integer $n\ge3$, let $L_n$ be the graph consisting of two paths $x_1x_2...x_n$, $y_1y_2...y_n$ and a matching $x_1y_1, ..., x_ny_n$. Let $A_n$ and $A_n'$ (see Figure \ref{fig:a}) be the two graphs obtained from $L_n$ by adding a pair of edges $x_1x_n, y_1y_n$, or $x_1y_n, y_1x_n$, respectively. Note that $A_3$ is the prism, $A_4$ is the cube, and $A_3'=K_{3,3}$. In addition, $A_n$ is the planar dual of $B_n$ while $A_n'$ is nonplanar. Let $\mathcal A=\{A_n,A_n':n\ge4\}$.

\begin{theorem}[\cite{w4c}] \label{thm:geelen}
For every weakly $4$-connected graph $G$ there exists a $(G,\{K_4,\Pi\}\cup\mathcal A\cup\mathcal B^\diamond\})$-chain of weakly $4$-connected graphs $G_0,G_1,...,G_t$ such that, for each $i=1,...,t$, either  $\|G_i\|=\|G_{i-1}\|-1$ or $G_i=G_{i-1}\del e/f$, where $G_{i-1}$ has a paw containing both $e$ and $f$.
\end{theorem}

In \cite{w4c}, this theorem was established and stated for general matroids. The formulation provided here is a straightforward restriction of this general result to graphic matroids. We will not get into a detailed matroidal discussion on how the restriction is obtained, since it would require some preparations on matroid theory. Instead, we provide a few remarks for readers who are interested in these details. First, it is well known that a connected graph is 3-connected if and only if its graphic matroid is 3-connected (see 8.1.9 of \cite{oxley}). Also, for any 3-connected graph $G$ and any $k\ge4$, the graphic matroid of $G$ has a 3-separation $(E_1,E_2)$ with $\min\{|E_1|,|E_2|\}\ge k$ if and only if $G$ has a 3-separation $(G_1,G_2)$ with $\min\{\|G_1\|,\|G_2\|\}\ge k$ (see 4.1 of \cite{bonds}). Therefore, a connected graph is weakly 4-connected if and only if its graphic matroid is weakly 4-connected. Moreover, a 4-element 3-separating set of the graphic matroid of $G$ is exactly a paw of $G$. Finally, as we remarked earlier, a graphic matroid is a trident if and only if the graph is $\Pi$.

\begin{figure}[htb]
\centerline{\includegraphics[scale=0.38]{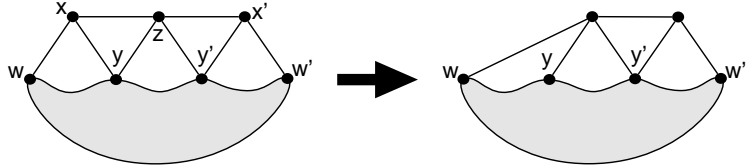}}
\caption{Operation $O_2$, where $d(x)=d(x')=3$ and $d(z)=4$.}
\label{fig:w4c}
\end{figure}

Let $K_6^-$ denote the graph obtained from $K_6$ by deleting an edge and let $K_{3,3}^+$ denote the graph obtained from $K_{3,3}$ by adding an edge. Let $\mathcal A_3=\mathcal A\cup\{A_3,A_3'\}$ and $\mathcal B_3=\mathcal B\cup \{B_3,B_3^+\}$. Our strengthened chain theorem requires the following operation. Let $G$ contain seven distinct vertices as shown in Figure \ref{fig:w4c}, where $d(x)=d(x')=3$ and $d(z)=4$. Let $O_2(G)$ denote $G/xz\del  xy$.

\begin{theorem} \label{thm:w4c}
For every weakly $4$-connected graph $G$ not contained in $\{K_4, W_4, K_6, K_6^-\}\cup\mathcal A_3\cup \mathcal B_3$ there exists a $(G, \{K_{3,3}^+,\Pi\})$-chain of weakly $4$-connected graphs $G_0,G_1,...,G_t$ such that, for each $i=1,...,t$, either $\|G_i\|=\|G_{i-1}\|-1$ or $G_i=O_2(G_{i-1})$.
\end{theorem}

We prove the theorem by proving a sequence of lemmas. We first prove that the 2-edge reduction stated in Theorem \ref{thm:geelen} can be replaced by the more restrictive operation $O_2$.

\begin{lemma}\label{lem:w4c}
Suppose a paw $P$ of a weakly $4$-connected graph $G$ contains two edges $e,f$ such that $G\del e/f$ is weakly $4$-connected. Then at least one of the following holds: \\
\indent (i) \ $G=W_4$; \\
\indent (ii) \ $G$ has a weakly $4$-connected minor with $\|G\|-1$ edges; \\
\indent (iii) \ $G\del e/f=O_2(G)$.
\end{lemma}

\noindent Proof. Since  $G\del e/f$ is weakly $4$-connected, we must have $|G\del e/f|\ge4$, which implies $|G|\ge 5$. Let $wx$, $xy$, $xz$, $yz$ be the four edges of $P$.

\noindent{\bf Claim 1.} {\it If  $yw$ or $zw$ is an edge of $G$ then either (i) or (ii) holds.}

By symmetry, suppose $yw\in G$. Let $G_1=P+yw$ and let $G_2$ be formed by all other edges of $G$. Then $(G_1, G_2)$ is a 3-separation of $G$. Since $G$ is weakly 4-connected and $\|G_1\|\geq 5$, we must have $\|G_2\|\le4$. Since very vertex outside $G_1$ is incident with at least three edges, any two of such vertices must be incident with at least five distinct edges. Therefore, $|G_2|=4$. Now it is clear that $G=W_4$ or $B_3$, which proves Claim 1.

In the following we assume $yw,zw\not\in G$. Since $d(w)\ge3$, we must have $|G|\ge6$.

\noindent{\bf Claim 2.} {\it If $d(y)=3$ or $d(z)=3$ then (ii) holds.}

By symmetry, we assume $d(y)=3$. Let $v$ be the third neighbor of $y$. Then $v\ne w$ and thus $\{v,w,z\}$ is a cut of $G$. Since $G$ is weakly $4$-connected, we must have $|G|=6$ and further, $G=A_3$. Now we see that (ii) holds, which proves Claim 2.

In the following we assume $d_G(y)\geq 4$ and $d_G(z)\ge4$.

\noindent{\bf Claim 3.} {\it $G\del yz$ is 3-connected.}

Suppose on the contrary that $G\del yz$ has a 2-separation $(G_1,G_2)$. If some $G_i$ contains both $y,z$ then $(G_i+yz,G_{3-i})$ would be a 2-separation of $G$, which is impossible. Hence, by symmetry, we may assume $y\in G_1\del G_2$ and $z\in G_2\del G_1$. It follows that $V(G_1\cap G_2)$ consists of $x$ and another vertex $v$. Again, by symmetry, we assume $xw\in G_1$. Then, as $d_G(x)=3$, we have $N_{G_2}(x)=\{z\}$. In addition, since $d_G(z)\geq 4$, we see that $|G_2|\geq 4$.  Now, $\{v,z\}$ is a 2-cut of $G$, a contradiction, which proves Claim 3.

If  $G\del yz$ is weakly 4-connected then (ii) holds. Hence, we assume that $G\del yz$ is not weakly 4-connected. By Claim 3, $G\del yz$ has a 3-separation $(G_1,G_2)$ with $\|G_i\|\ge 5$ for both $i=1,2$. Since $(G_1,G_2)$ cannot be extended into a 3-separation of $G$, we may assume without loss of generality that $y\in G_1\del G_2$ and $z\in G_2\del G_1$. Consequently, $x\in G_1\cap G_2$.

\noindent{\bf Claim 4.} {\it $w\not\in G_1\cap G_2$.}

Suppose $w\in G_1\cap G_2$. Since $zw\not \in G$, $yw\not \in G$, $d_G(z)\geq 4$, and $d_G(y)\geq 4$, we must have $|G_1|\geq 5$ and  $|G_2|\geq 5$. It follows that both $(G_1\cup P, G_2\del x)$ and $(G_1\del x, G_2\cup P)$ are 3-separations of $G$. Therefore, $|G_1\del x| = |G_2\del x|=4$. Now it is straightforward to verify that $G$ is the pyramid $\Pi$. This is impossible since $G\del e/f$ is not weakly 4-connected. This contradiction proves Claim 4.

By Claim 4, we may assume without loss of generality that $w\in G_1\del G_2$. Then $N_{G_2}(x)=\{z\}$. Since $\|G_2\|\ge5$, we have $\|G_2\del x\|\ge4$ and thus $|G_2\del x|\ge4$. It follows that $(G_1\cup P,G_2\del x)$ is a 3-separation of $G$. As a result, we must have $\|G_2\del x\|=4=|G_2\del x|=d_G(z)$. Thus $G$ is as shown in Figure \ref{fig:w4c}. From $d_G(y)\ge 4$ and $d_G(w), d_G(y'), d_G(w')\ge 3$ we see that $\|G\|\ge12$.

Now we need to determine the locations for $e$ and $f$. Note that $f\ne yz$ since $G/yz\del e$ has a vertex of degree 2 (no matter how $e$ is chosen), which should not happen in a weakly 4-connected graph. Then it is routine to see that there are three possibilities: $G/xz\del e \cong G/xw\del xy$, $G/xy\del e\cong G/xw\del xz$, and $G/xw\del yz$. Let $H$ denote the resulting graph. In the last two cases, let $H_1=H\del\{z,x'\}$ and let $H_2$ be the subgraph of $H$ formed by the five edges incident with $z$ or $x'$. Then $(H_1,H_2)$ is a 3-separation of $H$. Moreover, since $\|H\|=\|G\|-2\ge10$, we conclude that $\|H_1\|\ge5$ and thus $H$ is not weakly 4-connected. Therefore, the first case must happen and thus $H=O_2(G)$, so (iii) holds. \qed

\begin{lemma}\label{lem:Ae}
If $G\del e\in\cal A$ then $G/f$ is weakly 4-connected for some $f$ disjoint from $e$.
\end{lemma}

\noindent Proof. Let $A=G\del e$ and $n=|G|/2$. Then $A\in\{A_n, A_n'\}$ and $G=A+e$. Let $x_1,...,x_n, y_1,...,y_n$ be the vertices of $A$ and let their adjacency be specified as in the definition of $A_n$ and $A_n'$. In this proof the indices are taken modulo $n$. We first consider the case $n=4$.

\noindent{\bf Claim 1.} {\it The lemma holds if $n=4$. }

\begin{figure}[htb]
\centerline{\includegraphics[scale=0.4]{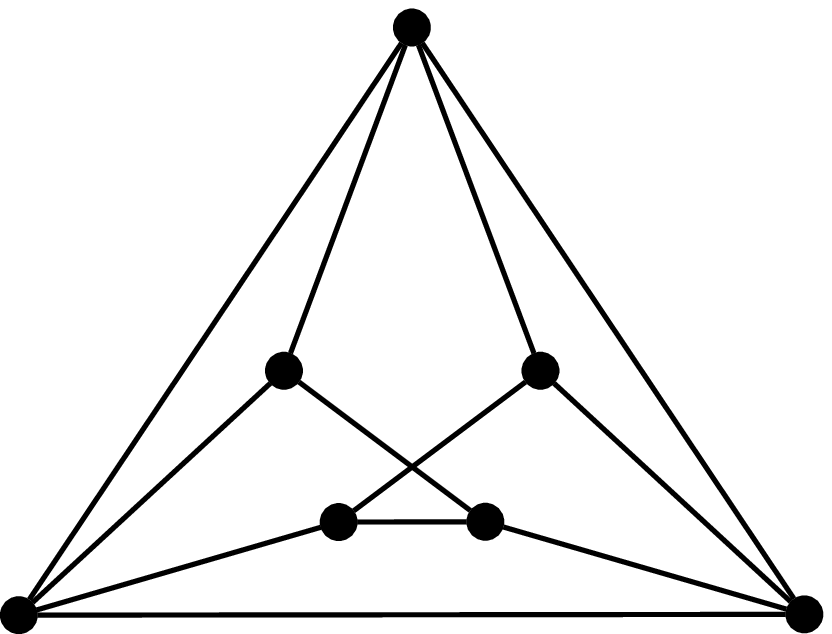}}
\caption{Kite $K$.}
\label{fig:kite}
\end{figure}

For $A=A_4$ or $A_4'$, it is straightforward to verify (separately) that, up to isomorphism, there are only two ways to add $e$ to $A$, and we may assume they are $y_3x_4$ and $y_3x_1$. So we have four cases, with two choices of $A$ and two choices of $e$. Let $f=x_2x_3$. It is routine to check that $(A_4+y_3x_4)/f\cong \Pi$ and in the other three cases $G/f\cong K$, which (known as a kite) is the graph shown in Figure \ref{fig:kite}. Since $\Pi$ and $K$ are weakly 4-connected, Claim 1 is proven.

In the following we assume $n\ge5$. To simplify our analysis we first ``place" $e$ in a special position.

\noindent{\bf Claim 2.} {\it We may assume without loss of generality that $e=y_3z$, where $z\not\in \{x_2,y_1\}$.}

Since $\{x_1y_1, ..., x_ny_n\}$ is a perfect matching of $A$, there must exist distinct $i$ and $j$ such that $e$ is between $\{x_i,y_i\}$ and $\{x_j,y_j\}$. Because of the rotational circular symmetry between $\{x_1,y_1\}, ...,  \{x_n,y_n\}$, we may assume $i=3$ and $j=i+t$, where $1\le t\le n/2$. Then $4\le j\le 3+\lfloor n/2\rfloor \le n$. Now, by the symmetry between $x_3$ and $y_3$, we may assume $e=y_3z$, where $z\in\{x_j,y_j\}$. Hence Claim 2 is satisfied.

In the following we prove that $f=x_2x_3$ satisfies the lemma. We first observe that $A/f$ is 3-connected, as deleting any vertex from it results in a 2-connected graph. Hence $G/f=A/f + e$ is 3-connected. Therefore, if $G/f$ is not weakly 4-connected, $G/f$ must have a 3-separation $(G_1,G_2)$ with $\|G_2\|\ge\|G_1\|\ge 5$. We need the following observation.

\noindent {\bf Claim 3.} {\it $G$ is weakly $4$-connected.}

For each vertex $v$ of $A$, note that $A\del v$ can be obtained from a 3-connected graph by subdividing each edge at most once, which implies that any 2-cut of $A\del v$ may cut off at most one vertex. It follows that every 3-cut of $A$ may cut off at most one vertex. That is, for every 3-separation $(H_1,H_2)$ of $A$, we must have $\min\{|H_1|,|H_2|\}=4$. Since $A$ has no triangles, we conclude that $\min\{\|H_1\|,\|H_2\|\}=3$, which implies $G=A+e$ is weakly 4-connected, as every 3-separation of $G$ induces a 3-separation of $A$.

\begin{figure}[htb]
\centerline{\includegraphics[scale=0.6]{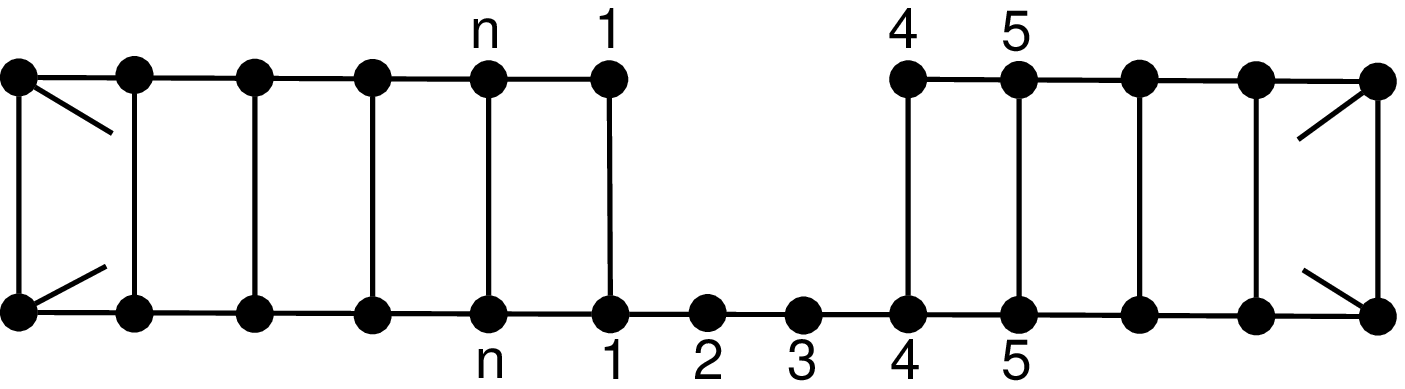}}
\caption{$A\del \{x_2, x_3\}$.}
\label{fig:q4cp}
\end{figure}

Claim 3 implies that, in the above 3-separation $(G_1,G_2)$, vertex $\bar f$ must be contained in both $G_1$ and $G_2$. Let $S$ denote the 3-cut $V(G_1\cap G_2)$. It follows that $S\del \bar f$ is a 2-cut of $G/f\del \bar f =G\del \{x_2,x_3\}=A\del \{x_2,x_3\} +e$. Since $n\ge5$, $A\del \{x_2,x_3\}$ is a subdivision of a 3-connected graph (see Figure \ref{fig:q4cp}). Thus any 2-cut of $A\del \{x_2,x_3\}$ can only separate some degree-2 vertices from the rest of the graph. It is straightforward to check that these cuts are: $\{y_1,y_3\}$, $\{y_1,y_4\}$, $\{y_2,y_4\}$, $\{y_4,x_5\}$, and $\{y_1,z_n\}$, where $z_n=x_n$ if $A=A_n$ and $z_n=y_n$ if $A=A_n'$. By Claim 2, $\{y_1,y_4\}$ and $\{y_2,y_4\}$ are not 2-cuts of $A\del \{x_2,x_3\} +e$. So we only need to consider the other three 2-cuts. If $S=\{\bar f, y_1,z_n\}$ or $S=\{\bar f, y_4,x_5\}$, then $G_1$ consists of the three edges incident with $x_1$ or $x_4$, respectively. If $S=\{\bar f, y_1,y_3\}$, then $G_1$ consists of the three edges  incident with $y_2$, and possibly also $y_3\bar f$ ($y_3y_1$ is not an edge, by Claim 2). In all cases, $\|G_1\|\ge5$ is not satisfied. This contradiction proves that $G/f$ is weakly 4-connected. \qed

The following can be considered as the dual of the last lemma.

\begin{lemma}\label{lem:q4csplit}
Let $e=xy$ be an edge of a $3$-connected graph $G$. If $G/e\in\cal B$ then $G$ has an edge $f$ such that $G\del f$ is weakly $4$-connected with $\delta(G\del f)=3$.
\end{lemma}

\noindent Proof. Let $B=G/e$. Then $B=B_n$ or $B_n^+$ for some $n\ge4$. Let $h_1$, $h_2$ be the two hubs of $B$ and let $C=v_1v_2...v_nv_1$ be the rim cycle. We first make two observations.

\noindent{\bf Claim 1.} {\it $B\del h_1v_1$ is quasi $4$-connected.}

Since $B$ is 4-connected, it follows that $B\del h_1v_1$ is 3-connected. Now let $(H_1,H_2)$ be a 3-separation of $B\del h_1v_1$. If $n=4$ then $|H_1|+|H_2|=|B|+3=9$, which implies $\min\{|H_1|,|H_2|\}\le4$. So we assume $n\ge5$. In this case it is routine to see that any two vertices different from $v_1$ are linked in $B\del h_1v_1\del h_2v_1$ by four independent paths. It follows that some $H_i$ has to contain all vertices different from $v_1$. Therefore, $H_{3-i}$ must contain precisely $v_1$ and its three neighbors, implying $\min\{|H_1|,|H_2|\}\le4$ again.

\noindent{\bf Claim 2.} {\it Let $uv$ be an edge of a $3$-connected graph $H$. If $H$ is not quasi $4$-connected but $H/uv$ is, then $H$ has a cubic vertex $w\ne u,v$ such that $|N_H(\{u,w\})|=3$ or  $|N_H(\{v,w\})|=3$. }

Since $H$ is 3-connected but not quasi 4-connected, it has a $3$-separation $(H_1,H_2)$ with $\min\{|H_1|, |H_2|\}\ge 5$. By symmetry, let $uv\in H_1$. Note that $|\{u,v\}\cap V(H_2)|\le 1$ because otherwise $H/uv$ would admit a 2-cut, which is not the case. It follows that $(H_1/uv,H_2)$ is a 3-separation of $H/uv$. Since $H/uv$ is quasi 4-connected and $d_{H/uv}(\overline{uv})\ge 4$, we must have $|H_1|=5$ and $|\{u,v\}\cap V(H_2)|=1$. Therefore, $V(H_1\del H_2)$ consists of a vertex $w$ and exactly one of $u$ and $v$. It follows that either $|N_H(\{u,w\})|=3$ or $|N_H(\{v,w\})|=3$. Finally, note that $3\le d_H(w)\le4$ since $H$ is 3-connected and $|H_1|=5$. If $d_H(w)=4$ then $uvw$ is a triangle, implying that $H/uv$ has a parallel edge, a contradiction. Thus, $d_H(w)=3$, which completes the proof of Claim 2.

\begin{figure}[htb]
\centerline{\includegraphics[scale=0.33]{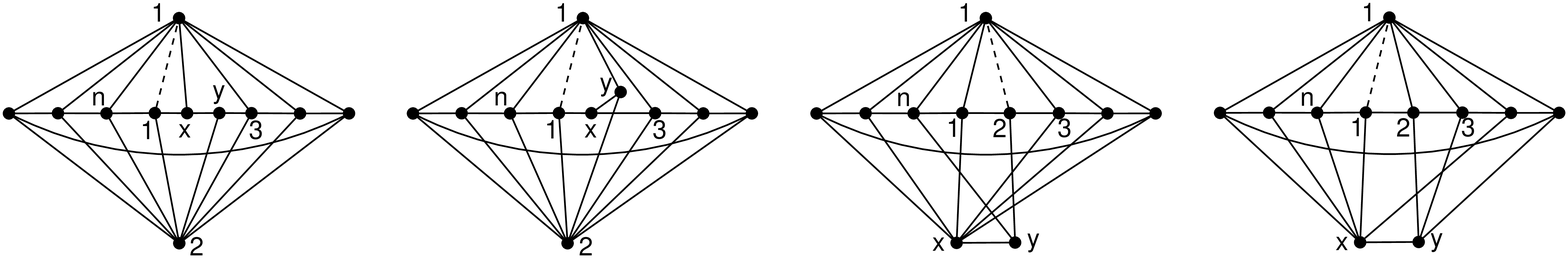}}
\caption{$H=G\del f$.}
\label{fig:w4cab}
\end{figure}

Now we prove the lemma by finding a required edge $f$. We first consider the case when $\bar e$ is a rim vertex. Without loss of generality, let $\bar e=v_2$. Up to isomorphism, there are two ways to split $v_2$, as shown in Figure \ref{fig:w4cab}. Let us assume by symmetry that $N_G(x)=\{y, v_1, h_1\}$ or $\{y, v_1, v_3\}$. We prove that $f=h_1v_1$ satisfies the lemma. Let $H=G\del f$. Note that each cubic vertex of $H$ (consisting of $x,y,v_1$, and possibly $h_1$) is contained in at most one triangle, so we only need to prove that $H$ is quasi 4-connected. By Claim 1, $H/e$ is quasi 4-connected. Then, by Lemma \ref{lem:split}, $H$ is 3-connected.
Hence, by Claim 2, if $H$ is not quasi 4-connected then $H$ has a cubic vertex $z\ne x,y$ such that $|N_H(\{t,z\})|=3$ for $t=x$ or $y$. Now it is routine to verify that no such $t,z$ exist, which proves that $H$ is weakly 4-connected.

It remains to consider the case when $\bar e$ is a hub of $B$. We assume that $B\ne B_4$ because in this case $\bar e$ can also be considered as a rim vertex and thus this is a case considered in the last paragraph. As a result, we have $d_B(h_i)\ge5$ ($i=1,2$). Without loss of generality, let us assume $\bar e=h_2$ and $d_G(x)\ge d_G(y)$. Then we have $d_G(x)\ge4$. We further assume, without loss of generality, that $xv_1,yv_2\in G$, as illustrated in Figure \ref{fig:w4cab}. If $N_G(y)=\{x,v_2,v_n\}$ then we choose $f=h_1v_2$, else we choose $f=h_1v_1$. We prove that our choice of $f$ satisfies the lemma. Let $H=G\del f$. If $f=h_1v_2$, neither of the two cubic vertices of $H$, $v_2$ and $y$, is contained in a triangle. If  $f=h_1v_1$, each cubic vertex of $H$ (consisting of $v_1$ and possibly $y$) is contained in at most one triangle. Therefore, in both cases, to establish the weak 4-connectivity of $H$ we only need to show that $H$ is quasi 4-connected. Let $i=1$ or 2 with $f=h_1v_i$. Note that $v_i$ is the only cubic vertex of $H$ outside $\{x,y\}$. By Claim 1, $H/e$ is quasi 4-connected. By Lemma \ref{lem:split}, $H$ is 3-connected. Hence, by Claim 2, if $H$ is not quasi 4-connected then $|N_H(\{v_i,x\})|=3$ or $|N_H(\{v_i,y\})|=3$. We need to deduce a contradiction.

Suppose $i=2$. Then we have $N_H(\{v_2,x\})=N_H(x)$ and $N_H(\{v_2,y\})=\{v_1,v_3,v_n,x\}$. It follows that $|N_H(\{v_2,x\})|=|N_H(x)|=d_G(x)\ge4$ and $|N_H(\{v_2,y\})|=4$, which give us the desired contradiction. Next, suppose $i=1$. Since $N_H(\{v_1,x\})\supseteq (N_H(x)\del v_1)\cup\{v_2\}$, we again have $|N_H(\{v_1,x\})|\ge |N_H(x)|=d_G(x)\ge4$. On the other hand, the choice of $f=h_1v_1$ implies $N_H(y)\ne N_H(v_1)$ and thus $N_H(\{v_1,y\})$ properly contains $N_H(v_1)$. Hence, $|N_H(\{v_1,y\})|>3$, and we obtain the desired contradiction again. \qed

\begin{lemma}\label{lem:w4ca}
Suppose $G$ is weakly $4$-connected and $G\not\in\{K_4,W_4,K_6^-,\Pi\}\cup\cal A\cup B$. Then $G$ must contain a weakly $4$-connected minor $H\not\in \{K_4\}\cup\cal A\cup B$ such that either $\|H\|=\|G\|-1$ or $H=O_2(G)$.
\end{lemma}

\noindent Proof. By Theorem \ref{thm:geelen}, $G$ has a weakly 4-connected minor $G'$ such that either $\|G'\|=\|G\|-1$ or $G'=G\del e/f$, where $G$ has a paw containing both $e,f$. Then, by Lemma \ref{lem:w4c}, we may assume that either $\|G'\|=\|G\|-1$ or $G'=O_2(G)$. It is clear that $G'\ne K_4$ since there is no 3-connected graph $H$ satisfying any of $H\del e=K_4$, $H/e=K_4$, or $O_2(H)=K_4$. If $G'\not\in\cal A\cup B$ then $H=G'$ satisfies the lemma. Hence, we assume $G'\in\cal A\cup B$ and we will use $G'$ to find $H$ that satisfies the lemma.

We first consider the case $G'\in\cal A$. If $G'=G/e$ or $O_2(G)$ then $G'$ would have a vertex of degree 4, which is not the case. So $G'=G\del e$. By Lemma \ref{lem:Ae}, $G$ has an edge $f$ disjoint from $e$ such that $G/f$ is weakly 4-connected. Since $G'$ is 3-regular, $G/f$ has three vertices ($\bar f$ and the two ends of $e$) of degree 4 and all its other vertices are cubic. It follows that $G/f\not\in\{K_4\}\cup \cal A\cup B$ and thus $H=G/f$ satisfies the lemma.

Next, assume $G'\in \cal B$. Note that $G'\ne O_2(G)$ because this would require $G'$ to have a cubic vertex, which is not the case. Hence, $G'=G\del e$ or $G/e$. We consider these two cases separately. Let $e=xy$ and let $h_1,h_2$ be the two hubs of $G'$.

Suppose $G'=G\del e$. It is clear $|G'|\ge6$. If $|G'|=6$ then $G'=B_4$ or $B_4^+$, which implies $G=B_4^+$ or $K_6^-$, respectively. This violates our assumption on $G$. Hence we must have $|G'|\ge7$. Since $G\not\in\cal B$, both $x,y$ are rim vertices of $G'$. Let $H=G\del xh_1$. Since $G'$ is 4-connected, $G$ is 4-connected and thus $H$ is 4-connected, as $H$ contains four independent paths between $x$ and $h_1$. Consequently, $H$ is weakly 4-connected. We also have $H\not\in\cal A$ since $H$ is not 3-regular. Furthermore, $H\not\in\cal B$ since $d_H(h_2)\ge 5= d_H(y)$ yet $H\del \{h_2,y\}$ is not a cycle. Therefore, $H$ satisfies the lemma.

Finally, let $G'=G/e\in\cal B$. By Lemma \ref{lem:q4csplit}, $G$ has an edge $f$ such that $G\del f$ is weakly 4-connected with $\delta(G\del f)=3$. Note that all vertices of $G\del f$, other than the ends of $e$ and $f$, must have degree exceeding 3. It follows that $G\del f\not\in\{K_4\}\cup \cal A\cup B$. Hence, $H=G\del f$ satisfies the lemma, which completes our proof. \qed

\noindent{\bf Proof of Theorem \ref{thm:w4c}}.
Let $G$ be weakly 4-connected and let $G \not \in \{K_4,W_4, K_6, K_6^-\}\cup\mathcal A_3\cup \mathcal B_3$. We assume $G\not\in \{K_{3,3}^+,\Pi\}$ because otherwise we are done. We only need to show that there exists $H$ satisfying\medskip\\
 ($*$) \ $H\in\{G\del e, G/e, O_2(G)\}$ (for some $e$), $H$ is weakly 4-connected, and $H\not\in \{K_4,W_4, K_6, K_6^-\}\cup\mathcal A_3\cup \mathcal B_3$\medskip\\
\noindent because repeatedly using this result will generate the desired chain.

By Lemma \ref{lem:w4ca}, $G$ has a weakly 4-connected minor $G'\not\in \{K_4\}\cup\mathcal A\cup \mathcal B$ such that either $\|G'\|=\|G\|-1$ or $G'=O_2(G)$. If $G'\not\in\{W_4, K_6, K_6^-,A_3,A_3',B_3,B_3^+\}$ then $H=G'$ satisfies ($*$). Hence, we assume that $G'$ is one of these seven. In the following we exam each of these seven graph and we either obtain a contradiction (indicating $G'$ cannot be that graph) or find $H$ satisfying ($*$).

If $G'=W_4$ then $G'=G\del e$ or $G/e$ ($G$ is too small for applying $O_2$). In the first case $G=B_3$ and in the second case $G=A_3$ or $A_3'$. Since both cases contradict $G\not\in\mathcal A_3\cup \mathcal B_3$, we conclude that $G'\ne W_4$.

If $G'=A_3$ or $A_3'$ then $G'=G\del e$, as $G'\in\{G/e,O_2(G)\}$ would imply that $G'$ has a vertex of degree exceeding 3, which is not the case. If $G'=A_3$, note that $G=A_3+e$ is not weakly 4-connected, as it has a cubic vertex belonging to two triangles. If $G'=A_3'$, note that $G=A_3'+e=K_{3,3}^+$, which violates the assumption $G\not\in \{K_{3,3}^+,\Pi\}$. These contradictions prove that $G'\not\in\{A_3,A_3'\}$.

If $G'=B_3$ then $G'=G\del e$ or $G/e$ ($G$ is too small for applying $O_2$). In the first case $G=B_3^+$ and in the second case $G=A_3+f$ (which is not weakly 4-connected) or $A_3'+f$ (which is $K_{3,3}^+$). These again contradicts our assumptions on $G$ and thus $G'\ne B_3$.

\begin{figure}[htb]
\centerline{\includegraphics[scale=0.5]{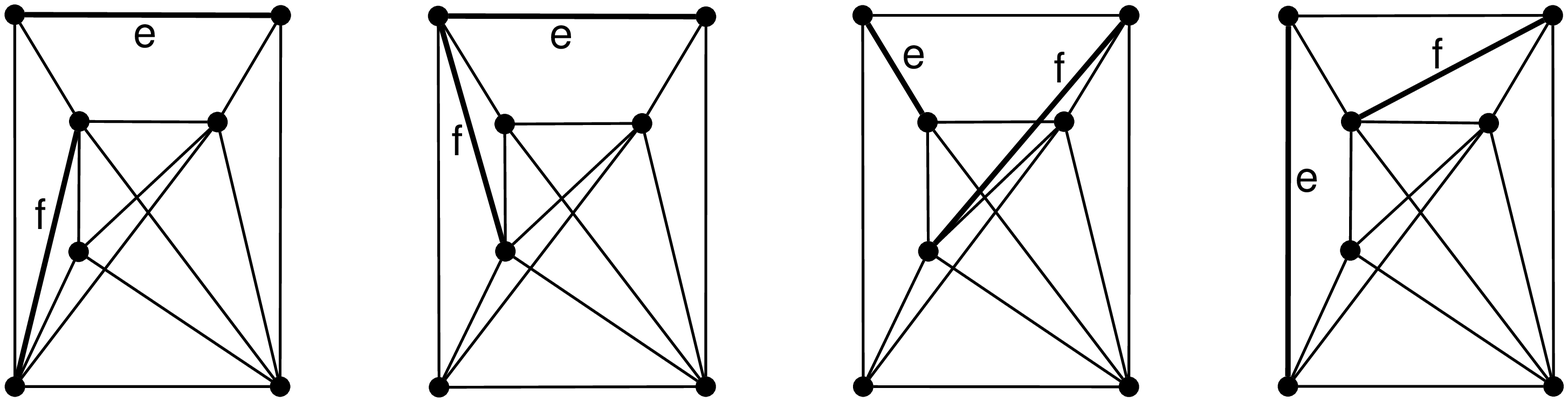}}
\caption{$G/e=K_6^-$ and $G\del f$ is weakly 4-connected.}
\label{fig:w4cabc}
\end{figure}

If $G'=B_3^+$ or $K_6$ then $G'=G/e$. To see this, note that $G'$ is a complete graph and thus $G'\ne G\del e$. Meanwhile, $G'\ne O_2(G)$ since $G'$ has no cubic vertices. Hence $G'=G/e$. If $G'=B_3^+$ then $G$ can also be obtained from $K_{3,3}^+$ by adding an edge $f$. In this case $H=G\del f$ satisfies ($*$).
If $G'=K_6$ then we may assume $N_G(v_1)=\{v_2, v_3, v_4\}$, $N_G(v_2)=\{v_1,v_5,v_6,v_7\}$, and $G\del\{v_1,v_2\}=K_5$. Then it is routine to verify that $H=G\del v_2v_7$ (isomorphic to the first graph in Figure \ref{fig:w4cabc}) satisfies ($*$).

Finally, suppose $G'=K_6^-$. Note that $G'\ne O_2(G)$ since $G'$ has no cubic vertices. If $G'=G\del e$ then $G=K_6$, contradicting our assumption on $G$. Hence $G'=G/e$. It is routine to verify that there are four choices for $G$, as illustrated in Figure \ref{fig:w4cabc}  (they can be obtained from $K_6$ by splitting a vertex and then deleting an appropriate edge). In all four cases, deleting $f$ results in the same graph, which we denote by $H$. Note that $H$ is weakly 4-connected (as $|N_H(\{u,v\})|\ge4$ for all $u,v$) and thus $H$ satisfies ($*$). This completes our proof of the theorem. \qed

\section{On quasi 4-connected graphs}

\begin{figure}[htb]
\centerline{\includegraphics[scale=0.46]{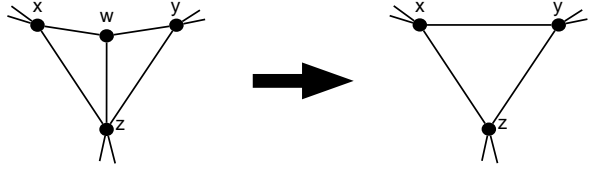}}
\caption{Operation $O_3$, where $d_G(x)\ge4$, $d_G(y)\ge4$, $d_G(z)\ge5$.}
\label{fig:q4c}
\end{figure}

In this section we obtain a similar result for quasi 4-connected graphs. We first state a known chain theorem. To do so, we need the following operation. Suppose $N(w)=\{x,y,z\}$, $xz,yz\in G$, $xy\not\in G$, $d(x),d(y)\ge4$ and $d(z)\ge5$, as shown in Figure \ref{fig:q4c}. Then $O_3(G)=G/xw\del wz$.

\begin{theorem}[\cite{q4c}] \label{thm:q4c}
For every quasi $4$-connected graph $G$ there exists a $(G,\{W_3,W_4,W_5\}\cup \cal A)$-chain of quasi $4$-connected graphs $G_0,G_1,...,G_t$ such that, for each $i=1,...,t$, either $\|G_i\| = \|G_{i-1}\|-1$ or $G_i=O_3(G_{i-1})$.
\end{theorem}

Let $\mathcal G_6$ denote the set of 3-connected graphs with at most six vertices. It is clear that every graph $G\in\mathcal G_6$ is quasi 4-connected since $G$ is too small to admit a 3-cut that separates two vertices of $G$ from another two vertices of $G$. Outside $\mathcal G_6$ the two smallest quasi 4-connected graphs are pyramid $\Pi$ and the graph $K$ shown below, which we call a {\it kite}. Both $\Pi$ and $K$ have seven vertices and twelve edges. The following is our strengthened chain theorem.

\begin{theorem}\label{thm:q4c+}
For every quasi $4$-connected graph $G\not\in \mathcal G_6\cup \mathcal A$ there exists a $(G$, $\{\Pi,K\})$-chain of quasi $4$-connected graphs $G_0$, $G_1$, ..., $G_t$ such that, for each $i=1,...,t$, either $\|G_i\| = \|G_{i-1}\|-1$ or $G_i=O_3(G_{i-1})$.
\end{theorem}

In the following we divide our proof into two lemmas. The first is a weaker version of the theorem. Recall that $A_3=$ prism is planar and $A_3'=K_{3,3}$ is nonplanar.

\begin{lemma}\label{lem:q4c+}
For every quasi $4$-connected graph $G$ not contained in $\{W_3$, $W_4$, $W_5$, $B_3$, $B_3^+\}\cup \mathcal A_3$ there exists a $(G$, $\{A_3,A_3'\})$-chain of quasi $4$-connected graphs $G_0$, $G_1$, ..., $G_t$ such that, for each $i=1,...,t$, either $\|G_i\| = \|G_{i-1}\|-1$ or $G_i=O_3(G_{i-1})$.
\end{lemma}

\noindent Proof. Let $G$ be quasi 4-connected and let $G\not\in \{W_3,W_4,W_5,B_3,B_3^+\}\cup \mathcal A_3$. We only need to show that there exists $H$ satisfying \medskip\\
\indent ($*$) \ $H\in\{G\del e, G/e, O_3(G)\}$ (for some $e$), $H$ is quasi 4-connected, and $H\not\in \{W_3, W_4, W_5, B_3, B_3^+\} \cup \mathcal A$ \medskip\\
because repeatedly using this result will generate the desired chain.

By Theorem \ref{thm:q4c}, $G$ has a quasi 4-connected minor $G'$ such that either  $\|G'\|=\|G\|-1$ or  $G'=O_3(G)$. We assume $G'\in \{W_3,W_4,W_5,B_3,B_3^+\}\cup \mathcal A$ because otherwise $H=G'$ satisfies ($*$). Depending on how $G'$ is obtained from $G$ we consider three cases. Each case will be further divided into subcases. In each subcase we will either obtain a contradiction (indicating that the subcase cannot happen) or find $H$ satisfying ($*$).

If $G'=O_3(G)$ then $G'$ contains at least three vertices of degree exceeding 3. It follows that $G' = B_3$ or $B_3^+$. By reversing operation $O_3$ we see that $G$ is a graph obtained from $K_{3,3}$ by adding, respectively, two edges $e,f$ or three edges $e,f,g$. Therefore, $H=G\del e$ satisfies ($*$).

If $G'=G/e$ then $d_{G'}(\bar e)>3$, which implies $G'\in\{W_4,W_5,B_3,B_3^+\}$. If $G'=W_4$ then $G=A_3$ or $A_3'$, contradicting our assumption on $G$. If $G'=B_3$ then $G$ can be expressed as $A_3+f$ or $A_3'+f$ for an edge $f$. In both cases $H=G\del f$ satisfies ($*$). If $G'=B_3^+$ then $G$ is a graph obtained from $K_{3,3}$ by adding two edges $f,g$. Consequently, $H=G\del f$ satisfies ($*$).
If $G'=W_5$, let $e=xy$ and let $v_1v_2v_3v_4v_5v_1$ be the rim cycle of $W_5$. Without loss of generality, let $N_{G}(x)=\{y,v_1,v_2,v_3\}$ or $\{y,v_1,v_3\}$. In both cases, $\{y,v_1,v_3\}$ is a 3-cut separating $\{x,v_2\}$ from $\{v_4,v_5\}$. It follows that $G$ is not quasi 4-connected, which is a contradiction.

It remains to consider the case $G'=G\del e$. Note that $G'\ne W_3$ or $B_3^+$ since they are complete graphs and thus they cannot be expressed as $G\del e$. If $G'=W_4$ or $B_3$ then $G=B_3$ or $B_3^+$, contradicting our assumption on $G$. If $G'\in\cal A$ then Lemma \ref{lem:Ae} ensures that $G$ has an edge $f$ such that $G/f$ is quasi 4-connected. Note that $G/f\not\in\cal A$ (since the degree of $\bar f$ is not 3) and $G/f \not \in \{W_3,W_4,W_5,B_3,B_3^+\}$ (since $|G/f|>6$). Therefore, $H=G/f$ satisfies ($*$). Finally, suppose $G'=W_5$. Let $v_1v_2v_3v_4v_5v_1$ be the rim cycle of $W_5$ and let $v_0$ be the hub. Without lose of generality, let $e=v_1v_3$. Now we see that $H=G\del v_0v_1$ (a graph obtained from prism $A_3$ by adding an edge) satisfies ($*$). The proof of the lemma is complete. \qed

Note that Theorem \ref{thm:q4c+} holds trivially for graphs with at most six vertices. In the following we prove that the theorem holds for graphs with seven vertices.

\begin{lemma}\label{lem:v6}
Let $G$ be quasi 4-connected with $|G|=7$. If $G$ is nonplanar then $G$ contains $K$ as a spanning subgraph; if $G$ is planar then $G$ contains $\Pi$ as a spanning subgraph.
\end{lemma}

\noindent Proof. We begin with two well known results on graph minors. (1) If a 3-connected graph is nonplanar then either the graph is isomorphic to $K_5$ or the graph contains a $K_{3,3}$ minor. (2) If a graph contains a cubic graph $Q$ as a minor then the graph contains a subdivision of $Q$ as a subgraph. From these two results and Lemma \ref{lem:q4c+} we immediately obtain: if $G$ is nonplanar then $G$ contains a subgraph $H$, where $H$ is a subdivision of $A=K_{3,3}$; if $G$ is planar then $G$ contains a subgraph $H$, where $H$ is a subdivision of $A=A_3$. We choose $H$ with $|H|$ maximized.

We first prove $V(H)=V(G)$. Suppose otherwise. Then $H\cong A$ (meaning that edges of $A$ are subdivided zero times) and $G$ contains exactly one vertex $z$ outside $H$. If two neighbors $x,y$ of $z$ are adjacent in $H$ then replacing edge $xy$ with path $xzy$ would result in a bigger $H$. Hence, $N_G(z)$ must be a stable set of $H$. Since $d_G(z)\ge3$, we must have $A=K_{3,3}$ and $N_G(z)$ has to be a color class of $H$. Now the quasi 4-connectivity of $G$ implies that $G\del N_G(z)$ contains a vertex $y$ of degree exceeding 1. It follows that $G\del y$ contains a subgraph $H'\cong A$ yet $N_G(y)$ is not stable in $H'$, which leads to a larger subdivided $A$ in $G$. This contradiction proves $V(H)=V(G)$.

\begin{figure}[htb]
\centerline{\includegraphics[scale=0.5]{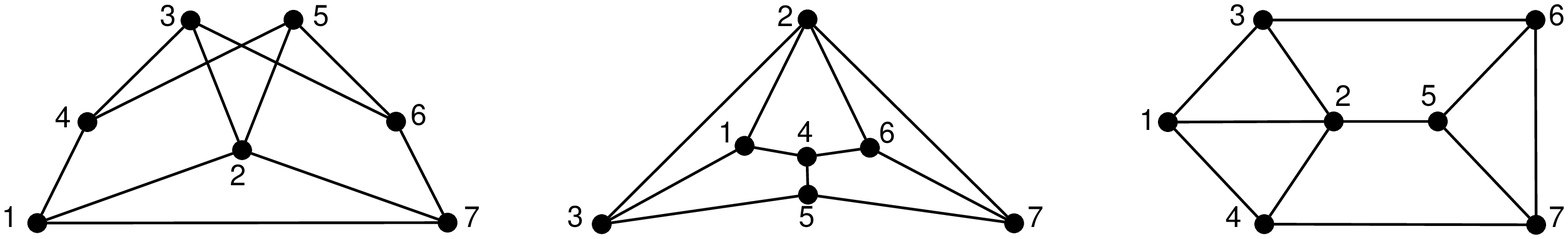}}
\caption{$H+v_1v_2$.}
\label{fig:v7}
\end{figure}

In the following discussion, let $V(G)=V(H)=\{v_1, ..., v_7\}$. Let $v_1$ be the unique subdividing vertex in obtaining $H$ and let $e=v_1v_2\in G$ be an edge outside $H$. Then $H+e$ is a subgraph of $G$. In the case $A=K_{3,3}$, up to isomorphism, $H+e$ is unique, which is the first graph in Figure \ref{fig:v7}. Since $\{v_2,v_4,v_6\}$ does not separate $v_3,v_5$ from $v_1,v_7$, $G$ must contain an edge $f$ between $\{v_3,v_5\}$ and $\{v_1,v_7\}$. By symmetry, there is only one choice for $f$ and it is easy to see that $H+e+f\cong K$, which settles the case for nonplanar $G$.

Next we assume $G$ is planar and $A=A_3$. If $v_1$ subdivides an edge of a triangle of $A$ then, subject to planarity and up to isomorphism, there are two choices of $e$. In these two cases we denote $H+e$ by $H_1$ and $H_2$, which are the last two graphs shown in Figure \ref{fig:v7}. If $v_1$ subdivides an edge not contained in a triangle of $A$ then, up to isomorphism, there is only one choice of $e$ and $H+e$ is isomorphic to $H_2$. In conclusion, $G$ contains $H_1$ or $H_2$ as a subgraph. If $H_1$ is a subgraph, since $\{v_2,v_4,v_5\}$ does not separate $v_1,v_3$ from $v_6,v_7$, $G$ must contain an edge $f$ between $v_1,v_3$ and $v_6,v_7$. By symmetry and planarity, there is only one choice for $f$ and it is easy to see that $H_1+f\cong \Pi$. Hence, we assume $H_1$ is not a subgraph, which implies that $H_2$ is a subgraph and (under the labeling of $H_2$) $v_1v_6,v_1v_7\not\in G$. Like in the previous cases, by considering cut $\{v_3,v_4,v_5\}$ we may assume $v_2v_6\in G$. Then by considering cut $\{v_2,v_4,v_6\}$ we have $v_3v_7\in G$, and by considering cut $\{v_2,v_3,v_7\}$ we also have $v_4v_5\in G$. Now we have $(H_2+v_2v_6+v_3v_7+v_4v_5)\del\{v_2v_3,v_6v_7\}\cong \Pi$, which completes our proof of the lemma. \qed

\medskip
\noindent{\bf Proof Theorem \ref{thm:q4c+}}. By Lemma \ref{lem:q4c+}, there exists a $(G,\{A_3, A_3'\})$-chain of quasi 4-connected graphs $G_0,G_1,...,G_t$ such that, fore all $i=1,...,t$, either $\|G_i\|= \|G_{i-1}\|-1$ or $G_i=O_3(G_{i-1})$. Since $|G_0|\ge7$, there exists a largest index $k$ with $|G_k|\ge 7$. Note that $k<t$, as $\|G_t\|=6$. By the maximality of $k$, we have $|G_{k+1}|<7$ and thus $G_{k+1}\ne G_k\del e$. Consequently, $G_{k+1}=G_k/e$ or $G_k/e\del f$ for some $e,f$. As a result, $|G_{k+1}| = |G_k|-1$, which implies $|G_k|=7$. By Lemma \ref{lem:v6}, $G_k$ contains either $\Pi$ or $K$ as a spanning subgraph. In other words, there exists a $(G_k,\{\Pi,K\})$-chain $H_0,...,H_s$ such that each $H_i$ ($i=1,...,s$) is obtained from $H_{i-1}$ by deleting an edge. Since $\Pi$ and $K$ are quasi 4-connected, every $H_i$ is also quasi 4-connected. Therefore, $G_0,...,G_k,H_1,...,H_s$ is a required chain, which proves the theorem. \qed

\section{On 3-connected graphs}

In this section we prove the following strengthening of Tutte Theorem. Since 3-connectivity has been extensively studied in the literature, it is possible to obtain a shorter proof for this result using more powerful tools. However, we choose to present a more elementary proof, which is not long either.

\begin{theorem} \label{thm:tutte+}
For every $3$-connected graph $G\not\in\cal W$ there exists a $(G,W_4)$-chain $G_0,G_1,...,G_t$ of $3$-connected graphs such that $||G_i||=||G_{i-1}||-1$ for all $i=1,...,t$.
\end{theorem}

\noindent
Proof. Let $G\not\in\cal W$ be 3-connected. We need to find an edge $f$ of $G$ such that at least one of $G\del f$ and $G/f$ is either the wheel $W_4$ or a 3-connected non-wheel graph. By Theorem \ref{thm:tutte}, we may assume that $G$ has an edge $e$ such that $G\del e$ or $G/e$ is a wheel $W_n$ with $n\ne 4$. Since no graph satisfies $G\del e=W_3$ or $G/e=W_3$, we may further assume $n\ge5$. Let $z$ be the hub of $W_n$ and let $e=xy$. If $G\del e = W_n$ then both $x,y$ are rim vertices. Let $f=xz$. Then $G\del f$ is not a wheel since it has two non-cubic vertices $y$ and $z$. In addition, $G\del f$ is 3-connected since $G$ is 3-connected and $G\del f$ has three independent paths joining $x$ and $z$. Therefore, $f$ satisfies the requirements. Next, suppose $G/e=W_n$. Then $G$ is obtained from $W_n$ by splitting $z$ into $x$ and $y$. Without loss of generality, let $d_G(x)\ge d_G(y)$. Then $d_G(x)>3$ since $n\ge5$. Let $f = uv$ be a rim edge such that $ux,vy\in G$. Then $f$ is not contained in any triangle and thus $G/f$ is simple. Now we see that $G/f$ is not a wheel since $\bar f$ and $x$ are two non-cubic vertices. In addition, by Lemma \ref{lem:split}, $H=G/f$ is 3-connected, as $si(H/e)=W_{n-1}$ is 3-connected. Therefore, $f$ satisfies the requirements. \qed

{\bf Acknowledgement:}   The second author is supported in part by NSF of China under grant 11961051 and supported in part by Natural Sciences Foundation of Guangxi Province under grant
2018GXNSFAA050117.

\end{document}